\documentclass[11pt,oneside]{amsart}

\usepackage[centertags]{amsmath}
\usepackage{amsfonts}
\usepackage{amssymb}
\usepackage{amsthm}
\usepackage{newlfont}

\newtheorem{thm}{Theorem}[section]
\newtheorem*{thm*}{Theorem}
\newtheorem{cor}[thm]{Corollary}

\newtheorem{lem}[thm]{Lemma}
\newtheorem{lemma}[thm]{Lemma}
\newtheorem{prop}[thm]{Proposition}
\newtheorem*{prop*}{Proposition}

\newtheorem*{conj*}{Conjecture}

\newtheorem*{dfn*}{Definition}
\theoremstyle{definition}
\newtheorem{rem}[thm]{\textbf{Remark}}
\newtheorem*{rmk*}{Remark}
\newtheorem*{fact*}{Fact}

\theoremstyle{proof}

\newcommand{\norm}[1]{\left\Vert#1\right\Vert}
\newcommand{\snorm}[1]{\Vert#1\Vert}
\newcommand{\abs}[1]{\left\vert#1\right\vert}
\newcommand{\set}[1]{\left\{#1\right\}}
\newcommand{\brac}[1]{\left(#1\right)}
\newcommand{\scalar}[1]{\left \langle #1 \right \rangle}

\newcommand{\Real}{\mathbb{R}}

\newcommand{\eps}{\varepsilon}

\newcommand{\Vol}[1]{\textnormal{Vol} \left(#1 \right)} 

\newcommand{\VolSet}[1]{\textnormal{Vol} \set{#1}}

\newcommand{\Var}{\textnormal{Var}}
\newcommand{\Prob}[1]{\text{Prob}\brac{#1}}

\newtheorem*{namedprop}{\theoremname}
\newcommand{\theoremname}{testing}

\newlength{\defbaselineskip}
\newcommand{\setlinespacing}[1]%
           {\setlength{\baselineskip}{#1 \defbaselineskip}}

\numberwithin{equation}{section}

\begin{document}

\title{On Gaussian marginals of uniformly convex bodies}

\author{Emanuel Milman}
\email{emanuel.milman@gmail.com}
\address{Department of Mathematics, The Weizmann Institute of Science, Rehovot 76100, Israel.}
\thanks{Supported in part by BSF and ISF}

\begin{abstract}
Recently, Bo'az Klartag showed that arbitrary convex bodies have
Gaussian marginals in most directions. We show that Klartag's
quantitative estimates may be improved for many uniformly convex
bodies. These include uniformly convex bodies with power type 2, and
power type $p>2$ with some additional type condition. In particular,
our results apply to all unit-balls of subspaces of quotients of
$L_p$ for $1<p<\infty$. The same is true when $L_p$ is replaced by
$S_p^m$, the $l_p$-Schatten class space. We also extend our results
to arbitrary uniformly convex bodies with power type $p$, for $2
\leq p < 4$. These results are obtained by putting the bodies in
(surprisingly) non-isotropic positions and by a new concentration of
volume observation for uniformly convex bodies.
\end{abstract}

\maketitle

\section{Introduction}

In recent years, numerous results have been obtained of the
following nature: let $X$ denote a uniformly distributed vector
inside a centrally-symmetric convex body $K$ of volume 1 in
$\Real^n$. Let $X_\theta := \scalar{X , \theta}$ denote its marginal
in the direction of $\theta \in S^{n-1}$, where $S^{n-1}$ denotes
the Euclidean unit sphere. Show that under suitable conditions on
$K$, the distribution of $X_\theta$ is \emph{approximately} Gaussian
for \emph{most} directions $\theta \in S^{n-1}$. Of course, the
meaning of ``approximately" and ``most" need to be carefully
defined, and vary among the different results.

To better illustrate this, consider the following examples. If $K =
[-\frac{1}{2},\frac{1}{2}]^n$, an $n$-dimensional cube, and $\theta
= \frac{1}{\sqrt{n}} (1,\ldots,1)$, the classical Central Limit
Theorem asserts that $\scalar{X , \theta}$ tends in distribution to
a Gaussian with variance $\frac{1}{12}$. Of course this is false for
\emph{all} directions $\theta \in S^{n-1}$, as witnessed by the
directions aligned with the cube's axes, but does hold for most
directions as measured by $\sigma$, the Haar probability measure on
$S^{n-1}$. When $K$ is a volume 1 homothetic copy of the Euclidean
ball $D_n$, the fact that (all) marginals are approximately Gaussian
is classical, dating back to Maxwell, Poincar\'{e} and Borel (see
\cite{DiaconisFreedmanDeFinetti} for a historical account). In the
broader context of general measures on $\Real^n$ with finite second
moment, Sudakov \cite{SudakovMarginals} showed that most marginals
are approximately the same \emph{mixture} of Gaussian distributions.
Under some additional conditions on the measure in question, Diaconis and Freedman
\cite{DiaconisFreedmanProjectionPursuit} showed that this mixture
can be replaced by a proper Gaussian. A generalized version of both
results was given by von Weizs{\"a}cker in \cite{vonWeiz}. Several
concrete convex bodies (other than the Euclidean Ball and the Cube),
such as the cross-polytope and simplex, were studied in
\cite{BrehmVoigt}.

Motivated by these and other results, it was conjectured by Antilla,
Ball and Perissinaki \cite{ABP} and Brehm and Voigt
\cite{BrehmVoigt} (using different and in fact stronger
formulations) that all convex bodies in $\Real^n$ have at least one
marginal which is approximately Gaussian, with the deviation tending
to 0 as the dimension $n$ tends to $\infty$.
This conjecture, referred to as the ``Central Limit Problem for
Convex Bodies" has been confirmed to hold for various classes of
convex bodies (\cite{ABP},\cite{Sodin-GaussianMarginals},
\cite{CoordCorrelationForOrliczNorms},\cite{Meckes-GaussianMarginals},\cite{KlartagEMilman-2-Convex}).

\medskip

Recall that $K$ is called \emph{isotropic} if it has volume 1 and
satisfies that $\Var(X_{\theta}) = L_K^2$ for all $\theta \in
S^{n-1}$ and some constant $L_K > 0$, which is called the isotropic
constant of $K$. Here $\Var(Y)$ denotes the variance of the random
variable $Y$. It is well known (e.g. \cite{Milman-Pajor-LK}) that every full-dimensional body has
an affine image which is isotropic and that this image is unique
modulo orthogonal rotations; we will refer to this affine image as
the body's isotropic position. Let us further denote the density
function of $X_\theta$ by $g_\theta(s) := \Vol{K \cap \set{s \theta
+ \theta^\perp}}$, and let $\phi_\rho(s) := \frac{1}{\sqrt{2 \pi}
\rho} \exp(-\frac{s^2}{2 \rho^2})$ denote the Gaussian density with
variance $\rho^2$. To emphasize that these notions depend on $K$, we
will usually use $g_\theta(K)$ instead of $g_\theta$, et cet.

\medskip

Recently, the Central Limit Problem for arbitrary convex bodies was
given an affirmative answer by Bo'az Klartag (\cite{KlartagCLP},
\cite{KlartagCLPpolynomial}) in the following sense: for every
isotropic convex body in $\Real^n$
\begin{equation} \label{eq:KlartagResults}
\sigma\set{ \theta \in S^{n-1} ; d_{TV}(g_\theta(K),\phi_{L_K}) \leq
\delta_n } \geq 1 - \mu_n ~,
\end{equation}
where $d_{TV}(f,g) = \int_{-\infty}^{\infty} \abs{f(s)-g(s)} ds$ is
the total-variation metric between the measures given by the
densities $f$,$g$, and $\delta_n$,$\mu_n$ are two series decreasing
to 0. Klartag's results in fact apply to all isotropic log-concave
probability measures on $\Real^n$. We refer to
\cite{Borell-logconcave} for the definition of log-concave measures,
and only remark that the Gaussian measure and arbitrary marginals of
convex bodies are known to be log-concave. In addition, for suitable
$k=k(n)$ increasing with $n$, the existence of $k$-dimensional
marginals which are approximately Gaussian was also shown. In
\cite{KlartagCLP} and later in \cite{FleuryGuedonPaourisCLP},
$\delta_n$ and $k(n)$ were shown to have logarithmic dependence in
$n$, and in \cite{KlartagCLPpolynomial} this was improved to
polynomial dependence: there exists some $\kappa_1,\kappa_2>0$ such
that the results are valid for $\delta_n = n^{-\kappa_1}$ and $k(n) =
n^{\kappa_2}$. In addition, it was shown in
\cite{KlartagCLPpolynomial} that one may replace in
(\ref{eq:KlartagResults}) the metric $d_{TV}$ with the stronger
notion of proximity $d_{Lin}^T$, to be defined in
(\ref{eq:strong-proximity}), with $T = L_K n^{\kappa_3}$ and $\delta_n
= n^{-\kappa_4}$ for some $\kappa_3,\kappa_4>0$. According to
\cite{KlartagCLPpolynomial} and some recent improvement in
\cite[Section 7]{KlartagUnconditionalVariance}, one may use
$\kappa_1 = 1/60, \kappa_2 = 1/15, \kappa_3 = 1/24, \kappa_4 = 1/24$
in the above statements.

\medskip

In this note, which is based on a previous version
\cite{EMilmanGaussianMarginalsArxiv} posted on the arXiv before the
announcement of Klartag's results, we will focus on showing the
existence of approximately Gaussian marginals in a strong sense for
a rather wide class of symmetric convex bodies. Although our results
do not apply to general convex bodies as in Klartag's work, we are
able to obtain better quantitative bounds on the deviation between
the body's marginal and the corresponding Gaussian distribution (the
$\delta_n$ in (\ref{eq:KlartagResults})). Earlier results in this
direction which have been most influential to our work include
\cite{ABP}, \cite{Sodin-GaussianMarginals} and
\cite{KlartagEMilman-2-Convex}; other references are given later on.
In those and previously mentioned results, approximately Gaussian
marginals are found by requiring from $K$ that its volume be highly
concentrated around a thin spherical shell of radius $\sqrt{n}
\rho$, for some $\rho>0$ and $\eps<1/2$:
\begin{equation} \label{eq:concentration}
\Prob{ \abs{\frac{|X|}{\sqrt{n}} - \rho} \geq \eps \rho } \leq
\eps ~.
\end{equation}
Usually, in order to obtain this type of volume concentration, the
body $K$ is put in isotropic position.
Following \cite{KlartagEMilman-2-Convex} but contrary to other
approaches, and perhaps surprisingly, we will see in this note that
it turns out to be more useful to put the body $K$ in some
non-isotropic \emph{position} (or affine image), for which we can
show (\ref{eq:concentration}). We will say that $K$ is $D$-sub-isotropic if $K$
satisfies that $\Var(X_\theta) \leq D \rho^2$ for all $\theta \in
S^{n-1}$, where $D \geq 1$ is some fixed universal constant.

\medskip

Let us denote the average density over all possible directions by
$g_{avg}(s) := \int_{S^{n-1}} g_\theta(s) d\sigma(\theta)$. Let
$\rho_\theta^2$ denote the variance of the distribution
corresponding to the density $g_\theta$, and set $\rho_{max} =
\max_{\theta \in S^{n-1}} \rho_\theta$ and $\rho_{avg} =
\int_{S^{n-1}} \rho_\theta d\sigma(\theta)$. We reserve the symbols
$C$,$C'$,$C_1$,$C_2$,$c$,$c_1$,$c_2$, etc., to indicate positive
universal constants, independent of all other parameters, whose
value may change from one appearance to the next.

There are usually two steps in showing the existence of
approximately Gaussian marginals: first, show that $g_{avg}$ is
close to $\phi_\rho$, and then show that most densities $g_\theta$
are close to $g_{avg}$. Again, the meaning of ``close to" and
``most" vary between the results. In \cite{ABP}, the proximity
between two even densities $f_1,f_2$ was interpreted in a rather
weak sense, by using the Kolmogorov metric (for even densities):
\begin{equation} \label{eq:weak-proximity}
d_{Kol}(f_1,f_2) := \sup_{t>0} \abs{\int_{-t}^t f_1(s) ds -
\int_{-t}^t f_2(s) ds},
\end{equation}
which does not capture the similarity in the tail behaviour of the
densities. Note that when comparing a one dimensional
log-concave density with a Gaussian one, it is known (see
\cite[Theorem 3.3]{BrehmVoigtVogt}) that $d_{TV}$ and $d_{Kol}$ are
equivalent in the sense that $d_{Kol} \leq d_{TV} \leq h(d_{Kol})$
for some function $h(t) = O((t \log(1/t))^{1/2})$. In fact, Klartag
obtains some of his results in
\cite{KlartagCLP,KlartagCLPpolynomial} using $d_{Kol}$ and
translates them to $d_{TV}$ using the above remark. Hence all the
results stated in this note for $d_{Kol}$ can be easily translated
to the total-variation metric.

We summarize the two steps from \cite{ABP} into a single statement.
In fact, our first observation in this note is that the argument of
\cite{ABP}, originally derived for an isotropic body, applies to a
body in arbitrary position, with some penalty accounting for the
deviation from isotropic position, as measured by:
\[
C_{iso}(K) := \rho_{max}(K)/\rho_{avg}(K).
\]
This more general statement, which was already used (without
proof) in \cite{KlartagEMilman-2-Convex}, reads as follows:

\begin{thm}[Generalized from \cite{ABP}] \label{thm:extended-ABP}
Assume that (\ref{eq:concentration}) holds for a centrally-symmetric
convex body $K$ in $\Real^n$. Then for any $\eps<\delta<c_2$:
\begin{equation} \label{eq:ABP-gaussians}
\sigma\set{ \theta \in S^{n-1} ; d_{Kol}(g_\theta(K),\phi_\rho) \leq
\delta } \geq 1- C_1 C_{iso}(K) \sqrt{n} \log{n} \exp\brac{- \frac{c_3 n \delta^2}{C_{iso}(K)^2}} ~.
\end{equation}
\end{thm}

Theorem \ref{thm:extended-ABP} is proved in Section
\ref{sec:GaussianMarginals-1}. We remark that it is easy to check
that $c_1 \rho_{avg} \leq \rho \leq c_2 \rho_{avg}$ (for some
universal constants $c_1,c_2>0$), whenever $\rho$ satisfies
(\ref{eq:concentration}), so we will sometimes use $\rho_{max}(K) /
\rho$ in place of the above definition of $C_{iso}(K)$.

\medskip

In \cite{Sodin-GaussianMarginals}, Sasha Sodin interpreted the
proximity between two even densities $f_1,f_2$ in a much stronger
sense, by measuring the following Linnik type quantity (see
\cite{LinnikBook}):
\begin{equation} \label{eq:strong-proximity}
d_{Lin}^T(f_1,f_2) := \sup_{0\leq s \leq T}
\abs{\frac{f_1(s)}{f_2(s)} - 1} ~,
\end{equation}
where $T$ may be as large as some power of $n$. Of course, this
stronger notion requires a stronger condition on the
concentration of volume inside $K$:
\begin{equation} \label{eq:Sodin-concentration}
\Prob{ \abs{ \frac{|X|}{\sqrt{n}} - \rho } \geq t \rho } \leq A
\exp( - B n^\nu t^\tau )~,
\end{equation}
for all $0\leq t \leq 1$ and some $A,B,\nu,\tau>0$. In that
case, we summarize the two steps in \cite{Sodin-GaussianMarginals}
into the following single statement. The following formulation,
which is not difficult to check, extends Sodin's result, originally
formulated for bodies in $D$-sub-isotropic position (with the dependence
on $D$ implicit in the constants), to arbitrary convex bodies (by explicitly stating the
dependence on $D$ via the parameter $C_{iso}(K)$).

\begin{thm}[\cite{Sodin-GaussianMarginals}] \label{thm:Sodin}
Let $K$ denote a centrally-symmetric convex body in $\Real^n$ and
assume that (\ref{eq:Sodin-concentration}) holds. Given $0 < \delta < c$
and $\mu > 0$, set:
\[
T = \rho \min\brac{ \brac{ \frac{c n C_{iso}(K)^{-2}
\delta^4}{\log n + \log \frac{1}{\delta} + \mu} }^{\frac{1}{6}} ,
(c(A,B,\nu,\tau) \delta)^{\gamma/\nu} n^\gamma} ~,
\]
where $\gamma := \nu / (2 \max(\tau,1))$ and
$c(A,B,\nu,\tau)$ explicitly depends on $A,B,\nu,\tau$.
Then:
\begin{equation} \label{eq:Sodin-gaussians}
\sigma \set{\theta \in S^{n-1} ; d_{Lin}^T(g_\theta(K),\phi_\rho)
\leq \delta} \geq 1 - \exp(-\mu).
\end{equation}
\end{thm}

The key step in Klartag's results from
\cite{KlartagCLPpolynomial} was the confirmation that
(\ref{eq:Sodin-concentration}) holds for arbitrary isotropic convex
bodies (and more generally, log-concave densities) with $\nu =
0.33$, $\tau = 3.33$, $\rho=L_K$ and universal constants $A,B>0$.
Plugging this into Theorem \ref{thm:Sodin}, we see that
(\ref{eq:Sodin-gaussians}) holds for arbitrary isotropic convex
bodies with $T = L_K n^{\kappa_3}$ and $\delta = n^{-\kappa_4}$, for
e.g. $\kappa_3 = 1/24, \kappa_4 = 1/24$, as mentioned earlier.

\medskip

Klartag's approach to the Central Limit Problem for convex bodies,
being completely general, cannot exploit any good properties which
certain classes of convex bodies posses. Consequently, certain
results for concrete classes which preceded Klartag's solution,
still give better quantitative bounds. These classes can be roughly
divided into two categories.

The first contains convex bodies possessing certain symmetries;
these include the $l_p^n$ unit-balls
(\cite{ABP},\cite{Sodin-GaussianMarginals}), more generally
arbitrary unit-balls of generalized Orlicz norms
(\cite{CoordCorrelationForOrliczNorms}), or other types of
symmetries (\cite{Meckes-GaussianMarginals},
\cite{Meckes-GaussianMarginalsNew}). In a recent progress in this
direction, Klartag has obtained in
\cite{KlartagUnconditionalVariance} a Berry-Esseen type result for
the marginals of an arbitrary convex body symmetric with respect to
reflections about coordinate hyperplanes.

The second category contains classes of uniformly convex bodies
under certain restrictions
(\cite{ABP},\cite{Sodin-GaussianMarginals},\cite{KlartagEMilman-2-Convex}).
With any centrally-symmetric convex $K \subset \Real^n$ we associate
a norm $\norm{\cdot}_K$ on $\Real^n$. The modulus of convexity of
$K$ is defined as the following function for $0<\eps\leq 2$:
\begin{equation}
\delta_K(\eps) = \inf \left \{ 1 - \left \| \frac{x+y}{2} \right
\|_K \ ; \  \|x\|_K, \|y\|_K \leq 1, \| x - y \|_K \geq \eps \right
\}. \label{eq:modulus}
\end{equation}
Note that $\delta_K$ is affine invariant, so it does not depend on
the position of $K$. A body $K$ is called uniformly convex if
$\delta_K(\eps) > 0$ for every $\eps > 0$. A body $K$ is called
``$p$-convex with constant $\alpha$'' (see, e.g. \cite[Chapter
1.e]{LT-Book-II}), if for all $0 < \eps \leq 2$,
\begin{equation}
\delta_K(\eps) \geq \alpha \eps^p ~. \label{eq:p-convex-defn}
\end{equation}
It is known that in such case $p$ cannot be smaller than 2.

The restriction imposed on $p$-convex bodies is usually via an upper
bound on the diameter of $K$ in isotropic position (\cite{ABP}) or
more generally in sub-isotropic position
(\cite{Sodin-GaussianMarginals}). For a $2$-convex body $K$ with
constant $\alpha$, this restriction on the diameter in isotropic
position was previously removed by Klartag and the author in
\cite{KlartagEMilman-2-Convex}. This was achieved by using Theorem
\ref{thm:extended-ABP}, which as remarked above, holds in an
arbitrary position. By placing $K$ in L\"{o}wner's minimal diameter
position, it was shown that $diam(K) \leq C n^{1-\lambda}
/ \lambda$, where $\lambda>0$ depends only on $\alpha$, enabling control of the
deviation term $C_{iso}(K)$. In order to apply Theorem \ref{thm:extended-ABP}, the required
concentration (\ref{eq:concentration}) was then deduced using a concentration result of
M. Gromov and V. Milman \cite{Gromov-Milman} for uniformly convex bodies (as in \cite{ABP},\cite{KlartagEMilman-2-Convex}). In order to compare the result from
\cite{KlartagEMilman-2-Convex} with the results in this note, we provide it below:

\begin{thm}[\cite{KlartagEMilman-2-Convex}] \label{thm:2-convex}
Let $K \subset \Real^n$ denote a 2-convex body with constant
$\alpha$ and volume 1. Assume in addition that it is in L\"{o}wner's
minimal diameter position, and denote $\rho = \int_K |x| dx /
\sqrt{n}$. Then (\ref{eq:concentration}) holds with:
\[
\eps = c_1 \sqrt{\log n} \: \alpha^{-1/2} \lambda^{-1}
n^{-\lambda},
\]
where $\lambda = \lambda(\alpha) > 0$ depends on $\alpha$ only.
In addition, for any $\eps < \delta < c_2$:
\[
\sigma\set{ \theta \in S^{n-1} ; d_{Kol}(g_\theta(K),\phi_\rho) \leq
\delta } \geq 1- \exp\brac{-c_3 \alpha \lambda^2 n^{2\lambda} \delta^2 }.
\]
\end{thm}

Our second observation in this note is that the same argument works for
arbitrary $p$-convex bodies ($p>2$) which have a small type-$s$
constant for large enough $s$ (see Section
\ref{sec:GaussianMarginals-2} for definitions). It is easy to show
that such bodies have small diameter in L\"{o}wner's position, and
so the usual application of the Gromov--Milman concentration
gives the desired result. As for the case $p=2$, the penalty term
$C_{iso}(K)$ needs to be handled in order to apply Theorem
\ref{thm:extended-ABP}. We postpone the formulation of our general result (Theorem \ref{thm:p-convex}) until Section
\ref{sec:GaussianMarginals-2}
and only state the following corollary, pertaining to the unit-balls of subspaces of quotients
of two useful classes of normed spaces for $1<p<\infty$: $L_p$, the class of $L_p$-integrable functions on $[0,1]$, and $S_p^m$, the Schatten class of $m$ by $m$ complex or real matrices, equipped with the norm $\norm{A} = (tr
(A A^*)^{p/2})^{1/p}$.

\begin{thm} \label{thm:L_p}
Let $K$ denote the unit-ball of an $n$-dimensional subspace of
quotient of $L_p$ or $S_p^m$ for $1<p<\infty$, and assume it has
volume 1. Assume in addition that it is in L\"{o}wner's minimal
diameter position, and denote $\rho = \int_K |x| dx / \sqrt{n}$.
Then (\ref{eq:concentration}) holds with:
\[
\eps = c_1 \sqrt{rq} (\log n)^{\frac{1}{\max(p,2)}}
n^{-\frac{1}{r}},
\]
where $r = \max(p,q)$ and $q=p^*=p/(p-1)$. In addition, for any $\eps
< \delta < c_2$:
\[
\sigma\set{ \theta \in S^{n-1} ; d_{Kol}(g_\theta(K),\phi_\rho) \leq \delta } \geq 1- n^{\frac{5}{2}}
\exp\brac{- \frac{c_3}{rq} n^{\frac{2}{r}} \delta^2} ~.
\]
\end{thm}

With our extended formulation of Theorem \ref{thm:Sodin} at hand, we
can also give analogous results to those of Theorems
\ref{thm:2-convex}, \ref{thm:L_p} (and \ref{thm:p-convex} from Section \ref{sec:GaussianMarginals-2})
using the stronger notion of proximity between densities
(\ref{eq:strong-proximity}). Indeed, for $p$-convex bodies as above,
the Gromov--Milman argument already implies the stronger
concentration assumption (\ref{eq:Sodin-concentration}), and the
penalty of $C_{iso}(K)$ appearing in Theorem \ref{thm:Sodin} is
handled exactly as for the former notion of proximity. We will only state the analogue
of Theorem \ref{thm:L_p}, the analogue of Theorem \ref{thm:2-convex} (and \ref{thm:p-convex})
is stated in Section \ref{sec:GaussianMarginals-2}.

\begin{thm} \label{thm:L_p-strong}
With the same assumptions and notations as in Theorem \ref{thm:L_p},
(\ref{eq:Sodin-concentration}) holds with:
\[
 \nu = \min(2/q,1) \;,\; \tau = \max(p,2) \;,\; A = 4 \;,\; B = q^{-2} (cp)^{-p/2} ~.
\]
In addition,
(\ref{eq:Sodin-gaussians}) holds for any $0<\delta<c$ and $\mu > 0$ with:
\[
T = \rho \min\brac{ \brac{ \frac{ c \delta^4 (rq)^{-1}}{\log n + \log \frac{1}{\delta} + \mu} }^{\frac{1}{6}} n^{\frac{1}{3r}}, c(p)
\delta^{\frac{1}{\max(p,2)}} n^{\frac{1}{2r}}}.
\]
\end{thm}

In Section \ref{sec:GaussianMarginals-3}, we take on a different
approach, which relies on the results of Bobkov and Ledoux from
\cite{BobkovLedoux}. Contrary to other methods, which need to
control the \emph{global} Lipschitz constant of the Euclidean norm
$|x|$ w.r.t. $\norm{\cdot}_K$, the results in \cite{BobkovLedoux}
enable us to average out the \emph{local} Lipschitz constant of
$|x|$ on $K$. Unfortunately, our estimate for this average enables
us to deduce a result for $p$-convex bodies only in the range $2\leq
p < 4$. Surprisingly, the position of $K$ which we use
to obtain our bounds is ``half'' way between
the isotropic and the minimal mean-width positions (see Theorem \ref{thm:p-convex-variance}).
We state the result only using the stronger notion of
proximity $d_{Lin}^T$, an analogous version using the weaker
$d_{Kol}$ metric may also be deduced.

\begin{thm} \label{thm:2-4-convex}
Let $K \subset \Real^n$ denote a $p$-convex body with constant
$\alpha$ for $2 \leq p < 4$, and assume it has volume 1. Assume in
addition that it is in the position given by Theorem
\ref{thm:p-convex-variance} below, and denote $\rho^2 = \int_K |x|^2
dx / n$. Then (\ref{eq:Sodin-concentration}) holds with
\[
\nu = \frac{3}{8} - \frac{1}{2q} \;,\; \tau = \frac{1}{2} \;,\; A=2 \;,\; B =
c \alpha^{\frac{1}{2p}} /
\min(f(p,\alpha),\log(1+n))^{\frac{1}{2}} ~,
\]
where $q = p^* = p/(p-1)$ and $f$ is some implicit function (given by Lemma
\ref{lem:general-positions}). In addition,
(\ref{eq:Sodin-gaussians}) holds for any $0<\delta<c$ and $\mu > 0$ with:
\[
T = \rho \min\brac{ \brac{ \frac{c \alpha^{\frac{1}{p}}
\min(f(p,\alpha),\log(1+n))^{-1} \delta^4}{\log n + \log
\frac{1}{\delta} + \mu} }^{\frac{1}{6}} n^{\frac{1}{6p}}, (c(p,\alpha)
\delta)^{\frac{1}{2}} n^{\frac{3}{16}-\frac{1}{4q}}}.
\]
\end{thm}

Note that for the range $2 \leq p < 4$, the latter Theorem holds
without \emph{any} assumptions on the diameter of the $p$-convex
body (or the type-constant of the corresponding space). Even for
$p=2$, this is an improvement over Theorem \ref{thm:2-convex} which
was proved in \cite{KlartagEMilman-2-Convex} and Theorem
\ref{thm:2-convex-strong}, since there an implicit function $\lambda
= \lambda(\alpha)$ appears in several expressions and in particular
in the exponent of $n$ (in Theorem \ref{thm:2-4-convex} we can
always replace $f$ by $\log(1+n)$).

As a corollary, we strengthen Theorem \ref{thm:L_p} for unit-balls
of subspaces of quotients of $L_p$ or $S_p^m$ with $1<p<
\frac{16}{13}$, since in this range, $r$ in Theorem \ref{thm:L_p}
exceeds the value of $\frac{16}{3}$. These bodies are known to be $2$-convex
with constant $\alpha = c(p-1)$ (see Lemma \ref{lem:sqlp}), so we
may apply Theorem \ref{thm:2-4-convex}.
\begin{cor} \label{cor:2-4-convex}
Let $K$ be the unit-ball of an $n$-dimensional subspace of quotient
of $L_p$ or $S_p^m$ for $1<p\leq \frac{16}{13}$, and assume it has
volume 1. Assume in addition that it is in the position given by
Theorem \ref{thm:p-convex-variance} below, and denote $\rho^2 =
\int_K |x|^2 dx / n$. Then (\ref{eq:Sodin-concentration}) holds with:
\[
\nu = \frac{1}{8} \;,\; \tau = \frac{1}{2} \;,\; A=2 \;,\; B = c
(p-1)^{\frac{1}{4}} / \log(1+n)^{\frac{1}{2}} ~.
\]
In addition, (\ref{eq:Sodin-gaussians}) holds for any $0<\delta<c$ and $\mu > 0$ with:
\[
T = \rho \min\brac{ \brac{ \frac{c \alpha^{\frac{1}{2}}
\log(1+n)^{-1} \delta^4}{\log n + \log \frac{1}{\delta} + \mu}
}^{\frac{1}{6}} n^{\frac{1}{12}}, (c(p) \delta)^{\frac{1}{2}} n^{\frac{1}{16}}}.
\]
\end{cor}

\medskip

\noindent \textbf{Acknowledgments.} I would like to sincerely thank
my supervisor Prof. Gideon Schechtman for interesting discussions. I
would like to thank Bo'az Klartag for getting me interested in
the Central Limit Problem for convex bodies. I would also like to
thank Profs. S. Bobkov and M. Ledoux for answering my questions. Final thanks
go out to the anonymous referee, whose suggestions helped improved the presentation
of this note.

\section{Gaussian Marginals in Arbitrary Position} \label{sec:GaussianMarginals-1}

We dedicate this section to the proof of Theorem
\ref{thm:extended-ABP}, which was already used in
\cite{KlartagEMilman-2-Convex} to deduce Theorem \ref{thm:2-convex},
and which will be used in the next section for proving Theorems
\ref{thm:p-convex} and \ref{thm:L_p}.

\begin{proof}[Proof of Theorem \ref{thm:extended-ABP}]
We follow the proof in \cite{ABP}, emphasizing the necessary
changes. Denote $G(t) = \int_{S^{n-1}} \int_{-t}^t g_\theta(s) ds
d\sigma(\theta)$ and $\Phi_\rho(t) = \int_{-t}^t \phi_\rho(s) ds$.
It was shown in \cite{ABP} that under the condition
(\ref{eq:weak-proximity}):
\begin{equation} \label{eq:extended-ABP-avg-diff}
\abs{ G(t) - \Phi_\rho(t) } \leq 4 \eps + \frac{c}{\sqrt{n}}
\end{equation}
for any $t>0$, and this is still valid for any position of $K$
since the isotropicity of $K$ was not used in the argument at
all. Another important observation from \cite{ABP}, which holds
regardless of position, is that for every $t>0$, $\int_{-t}^t
g_\theta(s) ds$ is a reciprocal of a norm. More precisely,
denoting:
\[
\norm{x}_{t} = \frac{|x|}{\int_{-t}^t g_{\frac{x}{|x|}}(s) ds},
\]
it was shown in \cite{ABP} that $\norm{\cdot}_t$ is a norm for
any $t>0$ and that:
\begin{equation} \label{eq:tmp1}
a_t(\frac{x}{|x|}) |x| \leq \norm{x}_t \leq b_t(\frac{x}{|x|})
|x|,
\end{equation}
where $a_t,b_t$ satisfy for $\theta \in S^{n-1}$:
\begin{equation} \label{eq:tmp2}
a_t(\theta) = c_1 \max(\frac{\rho_\theta}{t},1) \;,\; b_t(\theta)
= c_2 \max(\frac{\rho_\theta}{t},1).
\end{equation}
To conclude that given $t>0$, the individual marginals
$\int_{-t}^t g_\theta(s) ds$ are close to their spherical mean
$G(t)$ (which is already known to be close to $\Phi_\rho(t)$), the
authors of \cite{ABP} invoke a classical result on concentration
of Lipschitz functions around their mean: if $f:S^{n-1}
\rightarrow \Real$ is a $\lambda$-Lipschitz function then:
\begin{equation} \label{eq:Lip-concentration}
\sigma\set{\theta \in S^{n-1} ; f(\theta) - \int_{S^{n-1}} f(\xi)
d\sigma(\xi) \geq \delta } \leq \exp( - C n \delta^2 / \lambda^2
).
\end{equation}
To this end, an estimate on the Lipschitz constant of
$\int_{-t}^t g_\theta(s) ds$ is needed. Unfortunately, a
straightforward application of the argument in \cite{ABP} (as
reproduced below) yields a Lipschitz constant of $C
\frac{\rho_{max}}{\rho_{min}}$, where $\rho_{min} = \min_{\theta
\in S^{n-1}} \rho_\theta$, and this is not good enough for our
purposes. We therefore modify the argument a little. For
$0<\gamma<1$, let:
\[
A_\gamma = \set{\theta \in S^{n-1} ; \rho_\theta \geq (1-\gamma)
\rho_{avg}}.
\]
Since $\rho_\theta^2 = \int_K \scalar{x,\theta}^2 dx$, it is clear
that $\rho_\theta$ is a norm in $\theta$, and therefore its
Lipschitz constant is bounded above by $\rho_{max}$. Hence by
(\ref{eq:Lip-concentration}):
\[
\sigma(A_\gamma) \geq 1 - \exp\brac{-\frac{C n
\gamma^2}{C_{iso}(K)^2}}.
\]
This means that for most directions, we can actually use $(1-\gamma)
\rho_{avg}$ as a lower bound on $\rho_\theta$. Let $a^\gamma_t :=
c_1 \max((1-\gamma)\rho_{avg} / t , 1)$, and define the modified
norm $\norm{x}_t^\gamma := \max(\norm{x}_t , a^\gamma_t |x|)$. Note
that by (\ref{eq:tmp1}) and (\ref{eq:tmp2}), we did not alter the
norm on $\theta \in A_\gamma$, for which $\int_{-t}^t g_\theta(s) ds
= 1/\norm{\theta}_t^\gamma$. As in \cite{ABP}, we evaluate the
Lipschitz constant of the latter expression:
\[
\abs{\frac{1}{\norm{\theta_1}_t^\gamma} -
\frac{1}{\norm{\theta_2}_t^\gamma}} \leq \frac{\norm{\theta_1 -
\theta_2}_t^\gamma}{\norm{\theta_1}_t^\gamma
\norm{\theta_2}_t^\gamma} \leq
\frac{b_t(\frac{\theta_1-\theta_2}{|\theta_1-\theta_2|})}{(a_t^\gamma)^2}
|\theta_1-\theta_2| \leq C \frac{C_{iso}(K)}{(1-\gamma)}
|\theta_1-\theta_2|,
\]
regardless of the value of $t$. Denoting $G^\gamma(t) =
\int_{S^{n-1}} \frac{1}{\norm{\theta}_t^\gamma} d\sigma(\theta)$,
(\ref{eq:Lip-concentration}) implies that:
\[
\sigma \set{ \theta \in S^{n-1} ;
\abs{\frac{1}{\norm{\theta}_t^\gamma} - G^\gamma(t)} \geq \eta }
\leq 2 \exp\brac{ - \frac{C n \eta^2 (1-\gamma)^2}{C_{iso}(K)^2}
}.
\]
Since $\frac{1}{\norm{\theta}_t^\gamma}$ and $\int_{-t}^t
g_\theta(s) ds$ are both bounded from above by absolute constants
and differ only outside the set $A_\gamma$, we have:
\begin{equation} \label{eq:mean-diff}
\abs{G^\gamma(t)-G(t)} \leq C' \sigma\set{\theta \notin A_\gamma}
\leq C' \exp\brac{-\frac{C n \gamma^2}{C_{iso}(K)^2}}.
\end{equation}
We can now conclude as follows. Let $\delta > 0$ be given, and
assume that $\delta$ is not greater than some absolute constant
$c>0$, so that we may define $\gamma = C_0 \delta < 1/2$. The fact
that $\rho_\theta$ is a norm implies (e.g.
\cite{Milman-Schechtman-Book}) that $\rho_{max} \leq C \sqrt{n}
\rho_{avg}$, and therefore choosing $C_0$ above big enough, we
always have by (\ref{eq:mean-diff}), $\abs{G^\gamma(t)-G(t)} \leq
\delta / 2$. Hence:
\begin{eqnarray}
\nonumber & & \sigma \set{\theta \in S^{n-1} ; \abs{ \int_{-t}^t
g_\theta(s) ds - G(t) } \geq \delta \text{ or } \theta \notin A_\gamma} \\
\nonumber & \leq & \sigma\set{ \theta \notin A_\gamma} +
\sigma\set{\theta \in S^{n-1} ;
\abs{\frac{1}{\norm{\theta}_t^\gamma} - G^\gamma(t)} \geq \delta
- \abs{G^\gamma(t)-G(t)}} \\
\nonumber &\leq& \exp\brac{ - \frac{C n \gamma^2}{C_{iso}(K)^2}} +
2 \exp\brac{ - \frac{C n (\delta/2)^2 (1-\gamma)^2}{C_{iso}(K)^2}
} \leq 3 \exp\brac{-\frac{C n \delta^2}{C_{iso}(K)^2}}.
\end{eqnarray}

Together with (\ref{eq:extended-ABP-avg-diff}), and denoting
$H_\theta(t) = \abs{\int_{-t}^t g_\theta(s) ds - \int_{-t}^t
\phi_\rho(s) ds}$, we have for each $t>0$:
\[
\sigma \set{\theta \in S^{n-1} ; H_\theta(t) \geq \delta +
4\eps + \frac{c}{\sqrt{n}} \text{ or } \theta \notin A_\gamma
} \leq 3 \exp\brac{-\frac{C n \delta^2}{C_{iso}(K)^2}}.
\]

To pass from this estimate to one which holds for all $t>0$
simultaneously, we use the same argument as in \cite{ABP}, by
``pinning" down $H_\theta(t)$ at $C \sqrt{n} \log(n) C_{iso}(K)$
points evenly spread on the interval $[0,C' \max(\rho,\rho_{max})
\log(n)]$. Since by our choice of $\gamma$, for $\theta \in
A_\gamma$ we have $\rho_\theta \geq \rho_{avg}/2$, it is easy to
verify (as in \cite{ABP}) that the Lipschitz constant of
$H_\theta(t)$ w.r.t. $t$ is bounded above by $C / \min(\rho_{avg} ,
\rho)$ on $A_\gamma$. By the remark after Theorem
\ref{thm:extended-ABP}, we know that $\rho$ and $\rho_{avg}$ are
equivalent to within universal constants, so the latter Lipschitz
constant is bounded above by $C' / \rho_{avg}$. Since the distance
between two consecutive ``pinned" points is $C \rho_{avg} /
\sqrt{n}$, this ensures that $H_\theta(t)$ does not change by more
than $C'' / \sqrt{n}$ between consecutive points, and this
additional error is absorbed by the earlier error terms. There is no
need to control $H_\theta(t)$ for $ t\geq C' \max(\rho,\rho_{max})
\log(n)$, since both $\int_t^\infty \phi_\rho(s) ds$ (Gaussian
decay) and $\int_t^\infty g_\theta(s) ds$ (log-concavity of
$g_\theta$, see Lemma 4 in \cite{ABP}), are smaller than $C /
\sqrt{n}$ in that range, and this is again absorbed by the previous
error terms. This concludes the proof.
\end{proof}

\section{Concentration of Volume in Uniformly Convex Bodies with Good
Type} \label{sec:GaussianMarginals-2}

In this section, we extend and strengthen the results from
\cite{KlartagEMilman-2-Convex} to $p$-convex bodies with ``good"
type. Recall that the (Rademacher) type-$p$ constant of a Banach
space $(X,\norm{\cdot})$ (for $1 \leq p \leq 2$), denoted
$T_p(X)$, is the minimal $T>0$ for which:
\[
\brac{\mathbb{E} \snorm{\sum_{i=1}^m \eps_i x_i}^2}^{1/2} \leq T
\brac{\sum_{i=1}^m \norm{x_i}^p}^{1/p}
\]
for any $m \geq 1$ and any $x_1,\ldots,x_m \in X$, where
$\set{\eps_i}$ are i.i.d. random variables uniformly distributed on
$\set{-1,1}$ and $\mathbb{E}$ denotes expectation.

As explained in the Introduction, the existence of Gaussian
marginals may be deduced using Theorems \ref{thm:extended-ABP} or
\ref{thm:Sodin}, once we show that the volume inside $K$ is
concentrated around a thin spherical shell, in some controllable
position of $K$. A fundamental observation on the concentration of
volume inside uniformly convex bodies was given by Gromov and Milman
in \cite{Gromov-Milman} (see also \cite{ABV} for a simple proof and
\cite{EMilmanSodinIsoperimetryForULC} for an isoperimetric version).
It states that if $K$ is uniformly convex with modulus of convexity
$\delta_K$, and $T \subset K$ with $|T| \geq \frac{1}{2} |K|$, then
for any $\eps > 0$:
\begin{equation} \label{eq:gm}
\frac{\Vol{ (T + \eps K) \cap K }}{\Vol{K}} \geq 1 - 2 \exp(-2 n
\delta_K(\eps)) ~.
\end{equation}
It is easy to see that the latter is equivalent to the
concentration around their mean of functions on $K$ which are
Lipschitz w.r.t. $\norm{\cdot}_K$.

Despite this attractive property of uniformly convex bodies, it is
still a hard task to deduce concentration of volume around some
spherical shell. The difficulty lies in the fact that for a convex
body $K$, the function $|x|$ has a Lipschitz constant of $diam(K)$
w.r.t. $\norm{\cdot}_K$, and this may be too big to be of use. In
the next section, we describe an approach for which we will only
need to control the \emph{average} Lipschitz constant of $|x|$ on
$K$, thereby eliminating the need to control $diam(K)$. In this
section, as in \cite{KlartagEMilman-2-Convex}, we use (\ref{eq:gm})
in a direct manner, by putting $K$ in a position for which we have
control over $diam(K)$. This will be ensured by the type condition
on $K$.

\medskip

We will use the following lemma, which is easy to deduce from
(\ref{eq:gm}) and the discussion above (see e.g. \cite{ABP} or
\cite[Lemma 5.2]{KlartagEMilman-2-Convex}):

\begin{lem} \label{lem:p-convex-concentration}
Let $K \subset \Real^n$ be a $p$-convex body with constant $\alpha$
and of volume 1. Then for any 1-Lipschitz (w.r.t. $|\cdot|$)
function $f$ on $K$:
\[
\VolSet{ x \in K ; \abs{ f(x) - \int_K f(y) dy } \geq diam(K) t } \leq 4
\exp(- 2 c^p \alpha n t^p).
\]
\end{lem}

Denoting $\rho = \int_K |x| dx / \sqrt{n}$ and $R = diam(K) /
\sqrt{n}$, we deduce:
\begin{equation} \label{eq:gen-concentration}
\VolSet{ x\in K ; \abs{ \frac{|x|}{\sqrt{n}} - \rho} \geq R t }
\leq 4 \exp( - 2 c^p \alpha n t^p ).
\end{equation}
We see that in order to get some non-trivial concentration, we
need to ensure that $R \ll n^{1/p}$.
We will make use of the following lemma from
\cite{EMilman-DualMixedVolumes} (which appeared first in an
equivalent form in \cite{Davis-etal-Lemma}):

\begin{lemma} \label{lem:type-2-lemma}
Let $K \subset \Real^n$ be a centrally-symmetric convex body in
L\"{o}wner's minimal diameter position. Then:
\[
M_2(K) diam(K) \leq T_2(X_K),
\]
where $M_2(K) = \brac{\int_{S^{n-1}} \norm{\theta}_K^2
d\sigma(\theta)}^{\frac{1}{2}}$.
\end{lemma}

By Jensen's inequality and polar integration, it is immediate for
a body of volume 1 that $M_2(K) \geq C / \sqrt{n}$, hence in
L\"{o}wner's position $diam(K) \leq \sqrt{n} T_2(X_K)$. By the
results from \cite{TJ-Type-Cotype-fev-vectors}, it is enough to
evaluate the type-$2$ constant of an $n$-dimensional Banach space
on $n$ vectors, and from this it is easy see that $T_2(X_K) \leq C
n^{\frac{1}{s} - \frac{1}{2}} T_s(X_K)$. We conclude:

\begin{cor} \label{cor:type-s-diam}
Let $K \subset \Real^n$ be a centrally-symmetric convex body of
volume 1 in L\"{o}wner's minimal diameter position. Then for any
$1\leq s \leq 2$:
\[
diam(K) \leq C n^{\frac{1}{s}} T_s(X_K).
\]
\end{cor}

Combining this with (\ref{eq:gen-concentration}), we immediately
have:

\begin{prop} \label{prop:type-concentration}
Let $K \subset \Real^n$ be a $p$-convex body with constant $\alpha$
and of volume 1. Assume in addition that it is in L\"{o}wner's
minimal diameter position, and denote $\rho = \int_K |x| dx /
\sqrt{n}$. Then for any $1\leq s\leq 2$ we have:
\begin{equation} \label{eq:type-concentration}
\VolSet{ x\in K ; \abs{ \frac{|x|}{\sqrt{n}} - \rho} \geq t }
\leq 4 \exp\brac{ - 2 c^p \alpha n^{1+\frac{p}{2} - \frac{p}{s}}
\brac{\frac{t}{T_s(X_K)}}^p }.
\end{equation}
\end{prop}

In order to get a meaningful result, i.e. a positive power in the
exponent of $n$, we see that we need to have a bounded type-$s$
constant $T_s(X_K)$ for $s > \frac{2p}{p+2}$. It was shown in
\cite{KlartagEMilman-2-Convex} that for a $2$-convex body $K$ with
constant $\alpha$, this is always satisfied for some $s = s(\alpha)
> 1$. More precisely, using the same notations as in
\cite{KlartagEMilman-2-Convex}, it was shown that there exists a
$0<\lambda<1/2$ depending solely on $\alpha$, such that for $s =
\frac{1}{1-\lambda}$ we have $T_s(X_K) \leq 1/\lambda$. By Corollary
\ref{cor:type-s-diam}, this means that a $2$-convex body $K$ with
constant $\alpha$, having volume 1 and in L\"{o}wner's position,
always satisfies:
\begin{equation} \label{eq:diam-2-convex}
diam(K) \leq C n^{1-\lambda}/\lambda.
\end{equation}
Plugging this into (\ref{eq:type-concentration}), we see that for
such a body:
\begin{equation} \label{eq:2-convex-concentration}
\VolSet{ x\in K ; \abs{ \frac{|x|}{\sqrt{n}} - \rho} \geq t }
\leq 4 \exp\brac{ - c \alpha n^{2\lambda} \lambda^2 t^2 }.
\end{equation}
Since in any position (e.g. \cite{Milman-Pajor-LK}):
\begin{equation} \label{eq:rho-lower-bound}
\rho \geq c_1 L_K \geq c_2,
\end{equation}
we get exactly the spherical concentration condition
(\ref{eq:Sodin-concentration}) needed for applying Theorem
\ref{thm:Sodin}. It remains to evaluate $C_{iso}(K)$, appearing in
Theorem \ref{thm:Sodin}. We argue as in
\cite{KlartagEMilman-2-Convex}, that $\rho_{max}(K)$ may be
evaluated just by examining the radii of the circumscribing ball of
$K$ and the inscribed Euclidean ball of $\tilde{K}=T(K)$, where $T$
is a volume preserving linear transformation such that $\tilde{K}$
is isotropic. Indeed, it is clear that $\rho_{max} =
\norm{T^{-1}}_{op} L_K$, where $\norm{\cdot}_{op}$ denotes the
operator norm. And if $K \subset R D_n$ and $\tilde{K} \supset r
D_n$, where $D_n$ denotes the Euclidean unit ball, it is clear that
$\norm{T^{-1}}_{op} \leq R/r$. In order to evaluate the radius of
the inscribed ball of $\tilde{K}$, we recall the following result
from \cite{KlartagEMilman-2-Convex}:

\begin{lem}[\cite{KlartagEMilman-2-Convex}] \label{lem:2-convex-inscribed-ball}
Let $K \subset \Real^n$ denote a $2$-convex body with constant
$\alpha$ and volume 1. If $K$ is in isotropic position then:
\begin{equation} \label{eq:2-convex-isotropic-radii}
c \sqrt{\alpha n} L_K D_n \subset K,
\end{equation}
implying in particular that $L_K \leq C / \sqrt{\alpha}$.
\end{lem}

Using (\ref{eq:diam-2-convex}) and Lemma
\ref{lem:2-convex-inscribed-ball}, we deduce that $\rho_{max}(K)
\leq C n^{\frac{1}{2}-\lambda} \alpha^{-1/2} \lambda^{-1}$. Using
(\ref{eq:rho-lower-bound}) and the remark after Theorem
\ref{thm:extended-ABP}, we conclude that:
\[
C_{iso}(K) \leq C n^{\frac{1}{2}-\lambda} \alpha^{-1/2}
\lambda^{-1}.
\]

Plugging everything into Theorem \ref{thm:Sodin}, we deduce:

\begin{thm} \label{thm:2-convex-strong}
Let $K \subset \Real^n$ denote a 2-convex body with constant
$\alpha$ and volume 1. Assume in addition that $K$ is in
L\"{o}wner's minimal diameter position, and denote $\rho = \int_K
|x| dx / \sqrt{n}$. Then (\ref{eq:Sodin-concentration}) holds with:
\[
\nu = 2 \lambda \;,\; \tau = 2 \;,\; A = 4 \;,\; B = c \alpha
\lambda^2 ~,
\]
where $\lambda = \lambda(\alpha) > 0$ depends on
$\alpha$ only.
In addition, (\ref{eq:Sodin-gaussians}) holds for any $0<\delta<c$ and $\mu > 0$ with:
\[
T = \rho \min\brac{ \brac{ \frac{c \alpha \lambda^2
\delta^4}{\log n + \log \frac{1}{\delta} + \mu} }^{\frac{1}{6}} n^{\lambda/3},
(c(\alpha) \delta)^{\frac{1}{4}} n^{\lambda/2}}.
\]
\end{thm}

We remark that Theorem \ref{thm:2-convex} was deduced in \cite{KlartagEMilman-2-Convex} by
choosing $t = c \sqrt{\log(n)} n^{-\lambda} \lambda^{-1}$ in
(\ref{eq:2-convex-concentration}) and applying Theorem
\ref{thm:extended-ABP}.

\medskip

For $p>2$ the situation is different, because $\frac{2p}{p+2}
> 1$ and we cannot in general guarantee that given $p$ and
$\alpha$, $T_s(X_K)$ is bounded even for $s=\frac{2p}{p+2}$. We will
therefore need to additionally impose some requirement on $T_s(X_K)$
for $s>\frac{2p}{p+2}$. Once this is done, we deduce from
(\ref{eq:type-concentration}), as for the case $p=2$, the spherical
concentration condition (\ref{eq:Sodin-concentration}) needed for
applying Theorem \ref{thm:Sodin}. In order to control the term
$C_{iso}(K)$ in this case, we need to generalize Lemma
\ref{lem:2-convex-inscribed-ball} to the case of $p$-convex bodies.
It is a mere exercise to repeat the proof in
\cite{KlartagEMilman-2-Convex}, which gives:

\begin{lem} \label{lem:p-convex-inscribed-ball}
Let $K \subset \Real^n$ denote a $p$-convex body with constant
$\alpha$ and volume 1. If $K$ is in isotropic position then:
\begin{equation} \label{eq:p-convex-isotropic-radii}
c (\alpha n)^{\frac{1}{p}} L_K D_n \subset K,
\end{equation}
implying in particular that $L_K \leq C
n^{\frac{1}{2}-\frac{1}{p}} \alpha^{-\frac{1}{p}}$.
\end{lem}
Arguing as above, this gives together with Corollary
\ref{cor:type-s-diam}:
\[
C_{iso}(K) \leq C n^{\frac{1}{s}-\frac{1}{p}}
\alpha^{-\frac{1}{p}} T_s(X_K).
\]
Plugging this together with Proposition
\ref{prop:type-concentration} into Theorem \ref{thm:Sodin}, we deduce:

\begin{thm} \label{thm:p-convex-strong}
Let $K \subset \Real^n$ denote a $p$-convex body with constant
$\alpha$ and volume 1. Assume in addition that $K$ is in
L\"{o}wner's minimal diameter position, and denote $\rho = \int_K
|x| dx / \sqrt{n}$. Then (\ref{eq:Sodin-concentration}) holds for
any $\frac{2p}{p+2} < s \leq 2$ with:
\[
\nu = 1 + p/2 - p/s \;,\;
\tau = p \;,\; A = 4 \;,\; B = \alpha (c / T_s(X_K))^p ~.
\]
In addition, (\ref{eq:Sodin-gaussians}) holds for any $0<\delta<c$ and $\mu > 0$ with:
\[
T = \rho \min\brac{ \brac{ \frac{c
\alpha^{\frac{2}{p}} T_s(X_K)^{-2} \delta^4}{\log n + \log
\frac{1}{\delta} + \mu} }^{\frac{1}{6}} n^{\frac{1}{6}+\frac{1}{3p} -
\frac{1}{3s}}, (c(p,\alpha,s)
\delta)^{\frac{1}{2p}} n^{\frac{1}{4}+\frac{1}{2p} -
\frac{1}{2s}}}.
\]
\end{thm}

Choosing:
\[
t = \frac{\log(n)^{1/p}T_s(X_K)}{c \alpha^{1/p} n^{\frac{1}{2} + \frac{1}{p} -
\frac{1}{s} }},
\]
we deduce from (\ref{eq:type-concentration}) the spherical
concentration condition (\ref{eq:concentration}) needed for applying
Theorem \ref{thm:extended-ABP}, and conclude:

\begin{thm} \label{thm:p-convex}
Let $K \subset \Real^n$ denote a $p$-convex body with constant
$\alpha$ and volume 1. Assume in addition that $K$ is in L\"{o}wner's
minimal diameter position, and denote $\rho = \int_K |x| dx /
\sqrt{n}$. Then (\ref{eq:concentration}) holds for
any $\frac{2p}{p+2} < s \leq 2$ with:
\[
\eps_s = c_1 T_s(X_K) (\log n)^{\frac{1}{p}}
\alpha^{-\frac{1}{p}} n^{-(\frac{1}{2} + \frac{1}{p} - \frac{1}{s})}.
\]
In addition, for any $\eps_s < \delta < c_2$:
\[
\sigma\set{ \theta \in S^{n-1} ; d_{Kol}(g_\theta(K),\phi_\rho) \leq
\delta} \geq 1-
n^{\frac{5}{2}} \exp\brac{- \frac{c_3 n^{1+\frac{2}{p}-\frac{2}{s}}
\delta^2 \alpha^{\frac{2}{p}}}{T_s(X_K)^2}}.
\]
\end{thm}

\medskip

It remains to deduce Theorems \ref{thm:L_p} and
\ref{thm:L_p-strong} about unit-balls of subspaces of quotients of
$L_p$ and $S_p^m$ for $1<p<\infty$. With Theorems
\ref{thm:p-convex-strong} and \ref{thm:p-convex} at hand, we only need to evaluate
these bodies' $r$-convexity and type-$s$ constants, for
appropriately chosen $r$ and $s$. This is done in the following
(essentially standard) lemma:

\begin{lem} \label{lem:sqlp}
Let $K \subset \Real^n$ denote the unit-ball of a subspace of
quotient of $L_p$ or $S_p^m$, for $1<p<\infty$. Let $r = \max(p,2)$,
$s = \min(p,2)$ and $q=p^*$. Then:
\begin{enumerate}
\item
$K$ is $r$-convex with constant $\alpha(p) = C
\min(p-1,p^{-1}2^{-p})$.
\item
$T_s(X_K) \leq C \max(\sqrt{p},\sqrt{q})$.
\end{enumerate}
\end{lem}
\begin{proof}[Sketch of Proof]
We will sketch the proof of the $L_p$ case. The proof of the $S_p^m$
case is exactly the same, since by the results of
N. Tomczak-Jaegermann \cite{TJ-Schatten-Modulus}, these two classes have equivalent type,
cotype and modulus of convexity (up to universal constants), and our
proof of the $L_p$ case will only depend on estimates for these
parameters.

It is known (e.g. \cite[Chapter 1.e]{LT-Book-II}) that up to
universal constants, $L_p$ has the same modulus of convexity as
$l_p$, and that the latter space is $r$-convex with constant
$\alpha(p)$. By definition, this is passed on to any subspace of
$L_p$, and it is easy to see that the same holds for any quotient
space (by passing to the dual and using the modulus of smoothness,
see \cite[Lemma 3.4]{KlartagEMilman-2-Convex}). Item (1) is thus
shown.

To show item (2), first consider the case $p \geq 2$. Since $L_q$ is
$2$-convex with constant $q-1$, the dual $L_p$ is $2$-smooth (see
\cite[Chapter 1.e]{LT-Book-II} or \cite{KlartagEMilman-2-Convex})
with constant $\beta = c(q-1)^{-1} \leq C p$, and by the above
discussion, the same is true for $K$ as a unit-ball of a subspace of
quotient of $L_p$. It is standard (e.g. \cite[Lemma
4.3]{KlartagEMilman-2-Convex}) that this implies that $T_2(X_K) \leq
C \sqrt{\beta} \leq C' \sqrt{p}$. When $p<2$, we use a different
argument. Denote by $C_q(X)$ the cotype-$q$ constant of a Banach
space $X$ and by $\norm{Rad(X)}$ the norm of the Rademacher
projection on $L_2(X_K)$ (see e.g. \cite{Milman-Schechtman-Book} for
definitions). Assuming that $K$ is the unit-ball of a subspace $S$
of a quotient $Q$ of $L_p$, we have:
\[
T_p(X_K) = T_p(S) \leq T_p(Q) \leq C \norm{Rad(Q)} C_q(Q^*),
\]
where the first inequality is immediate since type passes to
subspaces, and the second one is known (e.g.
\cite{Milman-Schechtman-Book}). But by duality, $Q^*$ is a subspace
of $L_q$, and therefore inherits the cotype-$q$ constant of $L_q$,
which is a universal constant (e.g. \cite{Milman-Schechtman-Book}).
We conclude that $T_p(X_K) \leq C \norm{Rad(Q)}$. But again by
duality $\norm{Rad(Q)} = \norm{Rad(Q^*)} \leq \norm{Rad(L_q)}$,
since $Q^*$ is a subspace of $L_q$. We use the standard estimates
$\norm{Rad(L_q)} \leq T_2(L_q) \leq C \sqrt{q}$ (e.g.
\cite{KlartagEMilman-2-Convex}) to deduce that $T_p(X_K) \leq C
\sqrt{q}$. This concludes the proof.
\end{proof}

Plugging this lemma into Theorems
\ref{thm:p-convex-strong} and \ref{thm:p-convex}, Theorems \ref{thm:L_p} and
\ref{thm:L_p-strong} are deduced.

\section{Concentration of Volume in $p$-Convex Bodies for $p<4$}
\label{sec:GaussianMarginals-3}

Let $K$ denote a $p$-convex body in $\Real^n$. As already
mentioned, it was first noticed by Gromov and Milman
(\cite{Gromov-Milman}) that functions on $K$ which are Lipschitz
w.r.t. $\norm{\cdot}_K$ are in fact concentrated around their
mean. This phenomenon has since been further developed by many
authors (e.g.
\cite{SchechtmanZinn},\cite{SchmuckyUniformlyConvexBodies},\cite{ABV}).
A common property to all of these approaches is that the level of
concentration depends on the global Lipschitz constant of the
function in question, even if in most places the function has a
much smaller local Lipschitz constant. The starting point in the
following discussion is the interesting results of Bobkov and
Ledoux in \cite{BobkovLedoux}, which overcome the above mentioned
drawback.

Recall that the \emph{entropy} of a non-negative function $f$
w.r.t. a probability measure $\mu$, is defined as:
\[
Ent_\mu(f) := \int f \log(f) d\mu - \int f d\mu \log( \int f
d\mu) ~.
\]
The expectation and variance of $f$ w.r.t. $\mu$ are of-course:
\[
E_\mu(f) := \int f d\mu \; , \; \Var_\mu(f) :=
E_\mu((f-E_\mu(f))^2) ~.
\]
We will also use the following notation for $q>0$:
\[
\Var^q_\mu(f) := E_\mu(|f-E_\mu(f)|^q) ~.
\]
We will use $Ent_K(f)$, $\Var_K(f)$ etc. when the underlying
distribution $\mu$ is the uniform distribution on $K$. We also
denote by $\norm{\cdot}^*$ the dual norm to $\norm{\cdot}$, defined
as $\norm{x}^* = \sup\set{ \abs{\scalar{x,y}} ; \norm{y} \leq 1}$.
The following log-Sobolev type inequality was proved in
\cite[Proposition 5.4]{BobkovLedoux} (we correct here a small
misprint which appeared in the original formulation):

\begin{thm}[\cite{BobkovLedoux}] \label{thm:Bobkov-Ledoux}
Let $K$ be a $p$-convex body with constant $\alpha$ and volume 1,
and let $q = p^* = p/(p-1)$. Then for any smooth function $f$ on
$K$:
\begin{equation} \label{eq:Bobkov-Ledoux}
Ent_K(\abs{f}^q) \leq \frac{2^q}{\Gamma(\frac{n}{p}+1)^{q/n}}
\brac{\frac{q}{\alpha}}^{q-1} \int_K (\norm{\nabla f}_K^*)^q dx ~.
\end{equation}
\end{thm}

When $p=q=2$, it is classical that this log-Sobolev type inequality
implies a Poincar\'e type inequality. Indeed, by applying Theorem
\ref{thm:Bobkov-Ledoux} to $f = 1+\eps g$ and letting $\eps$ tend to
0, we immediately have:
\[
\Var_K(g) \leq \frac{C}{\alpha n} \int_K (\norm{\nabla g}_K^*)^2 dx ~.
\]
More generally, it was shown in \cite{BobkovZegarlinski} that for
any $q \leq 2$ and norm $\norm{\cdot}$, a $q$-log-Sobolev type
inequality:
\[
\forall f \;\;\; Ent_\mu(|f|^q) \leq C \int \norm{\nabla f}^q
d\mu ~,
\]
always implies a $q$-Poincar\'e type inequality:
\[
\forall f \;\;\; \Var^q_\mu(f) \leq C \frac{2^q}{\log 2} \int
\norm{\nabla f}^q d\mu ~.
\]
Although with this approach the additional term $\frac{2^q}{\log 2}$
may not be optimal (as in the classical $q=2$ case), universal
constants do not play a role in our discussion. Applying this
observation to the $q$-log-Sobolev inequality in Theorem
\ref{thm:Bobkov-Ledoux} we deduce:
\begin{cor} \label{cor:Bobkov-Ledoux-Poincare}
With the same notations as in Theorem \ref{thm:Bobkov-Ledoux}:
\begin{equation} \label{eq:Bobkov-Ledoux-Poincare}
\Var^q_K(f) \leq \frac{C}{(\alpha n)^{q-1}} \int_K (\norm{\nabla
f}_K^*)^q dx ~.
\end{equation}
\end{cor}

Our goal will be to show some non-trivial concentration of the
function $g = |x|^2$ around its mean, which is tantamount to the
concentration of volume inside $K$ around a thin spherical shell. The advantage
of the estimates in Theorem \ref{thm:Bobkov-Ledoux} and Corollary
\ref{cor:Bobkov-Ledoux-Poincare} is that they ``average out" the
local Lipschitz constant of $f$ (w.r.t. $\norm{\cdot}_K$) at $x \in
K$, which is precisely $\norm{\nabla f(x)}_K^*$. The usual way to
deduce exponential concentration of $g$ around its mean is via the
Herbst argument, by applying Theorem \ref{thm:Bobkov-Ledoux} to the
function $f = \exp(\lambda g/q)$ (see \cite{BobkovLedoux} or
\cite{BobkovZegarlinski}) and optimizing over $\lambda$.
Unfortunately, estimating the right-hand side of
(\ref{eq:Bobkov-Ledoux}) for the function $\exp(\lambda |x|^2/q)$ is
a difficult task. An alternative way, which will a-priori only produces polynomial
concentration of $g$ around its mean, is to apply Corollary
\ref{cor:Bobkov-Ledoux-Poincare} to the function $f=g$ and use
Markov's inequality, in hope that estimating the right-hand side of
(\ref{eq:Bobkov-Ledoux-Poincare}) should be easier for $g$ itself. We will see that this
will in fact lead to exponential bounds.
We remark that it is possible to do the same with $f=g$ in
(\ref{eq:Bobkov-Ledoux}) and gain an additional logarithmic factor
in the resulting concentration, but we avoid this for simplicity. We
therefore start by applying Corollary
\ref{cor:Bobkov-Ledoux-Poincare} to the function $f = |x|^2$:

\begin{equation} \label{eq:Poincare-for-norm}
\Var^q_K(|x|^2) \leq \frac{C'}{(\alpha n)^{q-1}} \int_K
(\norm{x}_K^*)^q dx ~.
\end{equation}

In the following Proposition we estimate the right-hand side of
(\ref{eq:Poincare-for-norm}). We denote by $M^*(K)$ half the
mean-width of $K$, i.e. $M^*(K) = \int_{S^{n-1}} \norm{\theta}_K^*
d\sigma(\theta)$. We also denote by $SL(n)$ the group of volume
preserving linear transformations in $\Real^n$.

\begin{prop} \label{prop:variance-concentration}
Let $K$ be a $p$-convex body with constant $\alpha$. Assume that
$K$ is isotropic and of volume 1, and set $q=p^*=p/(p-1)$. Then
for any $T \in SL(n)$:
\[
(\Var^q_{T(K)}(|x|^2))^{\frac{1}{q}} \leq \frac{C'}{(\alpha
n)^{\frac{1}{p}}} n^{3/4} M^*(T^* T (K) ) L_K ~.
\]
\end{prop}
\begin{proof}
Since:
\[
\int_{T(K)} (\norm{x}_{T(K)}^*)^q dx = \int_{K} (\norm{x}_{T^*
T(K)}^*)^q dx ~,
\]
by (\ref{eq:Poincare-for-norm}) and a standard Lemma of C. Borell
\cite{Borell-logconcave} (note that $q \leq 2$):
\[
(\Var^q_{T(K)}(|x|^2))^{1/q} \leq \frac{C'}{(\alpha
n)^{\frac{1}{p}}} (\int_{K} (\norm{x}_{T^* T(K)}^*)^q dx)^{1/q}
\leq\frac{C''}{(\alpha n)^{\frac{1}{p}}} \int_{K} \norm{x}_{T^*
T(K)}^* dx ~.
\]
Let us evaluate the integral on the right. First, notice that the
contribution of $\set{x \in K \setminus C \sqrt{n} L_K D_n}$ to this
integral is negligible. To show this, we turn for simplicity to a
recent result of Grigoris Paouris (\cite{Paouris-IsotropicTail}),
who showed that when $K$ is in isotropic position:
\[
\Vol{ K \setminus C \sqrt{n} L_K t D_n} \leq \exp( - \sqrt{n} t )
\]
for all $t \geq 1$, hence:
\[
\int_{K \setminus C \sqrt{n} L_K D_n} \norm{x}_{T^* T(K)}^* dx
\leq \exp( - \sqrt{n}) diam(T^* T(K)) diam(K) ~.
\]
Since $diam(T^*T(K)) \leq C_1 \sqrt{n} M^*(T^*T(K))$ and $diam(K)
\leq C_2 n L_K$, we see that the latter integral is bounded by
$M^*(T^*T(K)) \exp( - \sqrt{n}/2 )$, which will be absorbed by the
estimate on the integral inside $K \cap C \sqrt{n} L_K D_n$. We
emphasize that neither the isotropic position nor Paouris'
estimate are cardinal here; a similar argument using Borell's
standard $\Psi_1$-estimate will give a negligible term. Denoting
$K' = T^*T(K)$, it remains to evaluate:
\begin{equation} \label{eq:integral-inside}
\int_{K \cap C \sqrt{n} L_K D_n} \norm{x}_{K'}^* dx ~.
\end{equation}
To this end, we apply a result of J. Bourgain (\cite{Bourgain-LK})
which uses the celebrated ``Majorizing-Measures Theorem" of
Fernique-Talagrand (see \cite{Talagrand-Book}), to deduce that the
latter is bounded by $C' n^{3/4} M^*(K') L_K$. We remark that this
is essentially the same argument which yields Bourgain's well known
bound on the isotropic constant $L_K \leq C n^{1/4} \log(1+n)$. For
completeness, we outline Bourgain's argument. The idea is to write
$\norm{x}_{K'}^*$ as $\sup_{y \in K'} \scalar{y,x}$, so
(\ref{eq:integral-inside}) becomes an expectation on a supremum of a
sub-Gaussian process. Let $X_H$ denote a random vector on the probability space $\Omega_H$ which is uniformly distributed on $K \cap C \sqrt{n} L_K D_n$, and for $y \in \Real^n$ denote $H_y := \scalar{X_H,y}$.
For a real-valued random variable $H$ on a probability space $(\Omega,d\omega)$
and $\alpha>0$, let $\norm{H}_{L_{\Psi_\alpha}(\Omega)}$ be defined
as:
\[
\norm{H}_{L_{\Psi_\alpha}(\Omega)} = \inf \set{ \lambda > 0 ;
\int_{\Omega} \exp ((H(\omega)/\lambda)^\alpha) d\omega \leq 2} ~.
\]
A standard calculation shows that:
\begin{eqnarray}
\nonumber \norm{H_y}_{L_{\Psi_2}(\Omega_H)} \leq
\sqrt{\norm{H_y}_{L_{\Psi_1}(\Omega_H)}
\norm{H_y}_{L_{\infty}(\Omega_H)}} \\
\nonumber \leq C_1 \sqrt{\norm{H_y}_{L_{2}(\Omega_H)}
\norm{H_y}_{L_{\infty}(\Omega_H)}} \leq C_2 n^{1/4} L_K |y| ~.
\end{eqnarray}
Denoting $H'_y = H_y / (C_2 n^{1/4} L_K)$, the latter implies that
the process $\set{H'_y}$ is sub-Gaussian w.r.t. the
Euclidean-metric, and hence by the Majorizing-Measures Theorem:
\[
E_{\Omega_H} \sup_{y \in K'} H'_y \leq C E_{\Omega_G} \sup_{y \in
K'} G_y ~,
\]
where $G_y := \scalar{X_G,y}$ and $X_G$ is a random vector on the probability
space $\Omega_G$ whose distribution is that of a standard $n$-dimensional Gaussian.
This implies that:
\[
E_{\Omega_H} \sup_{y \in K'} H_y \leq C_3 n^{1/4} L_K n^{1/2}
M^*(K') ~,
\]
and a similar bound holds for (\ref{eq:integral-inside}), since the
volume of $K \cap C \sqrt{n} L_K D_n$ is close to 1.
\end{proof}

It is easy to check (e.g. \cite{Milman-Pajor-LK}) that
$E_{T(K)}(|x|^2) = \int_{T(K)} |x|^2 dx \geq n L_K^2$, and therefore
any time the bound in Proposition \ref{prop:variance-concentration}
is asymptotically smaller than $n L_K^2$ we can deduce a
concentration result for $|x|^2$ on $K$. Unfortunately, we are
unable to do so in the isotropic position, which is perhaps the most
natural position for such concentration of volume to occur. For
example, when $K$ is a 2-convex isotropic body (with constant
$\alpha$), we cannot say much about $M^*(K)$; to the best of our
knowledge, the best upper bound was given in
\cite{KlartagEMilman-2-Convex}, where it was shown that in isotropic
position $M^*(K) \leq C(\alpha) n^{3/4}$, which is exactly the
critical value we wish to be properly below.

Proposition \ref{prop:variance-concentration} was deliberately
formulated in a way which enables us to work around this problem. We
will use a $T \in SL(n)$ so that $M^*(T^*T(K))$ is minimal. In order
to use Theorem \ref{thm:extended-ABP}, we will also need to control
$C_{iso}(T(K))$, which amounts (as in the previous section) to
controlling $\norm{T}_{op}$.
For 2-convex bodies, the relations between the isotropic, the
John and the minimal mean-width positions, were studied in
\cite{KlartagEMilman-2-Convex}. Recall that the John position of
a convex body $K$ is defined as the (unique modulo orthogonal
rotations) position with maximal radius of the inscribed
Euclidean ball. We summarize the additional relevant results from
\cite{KlartagEMilman-2-Convex} in the following:

\begin{lem}[\cite{KlartagEMilman-2-Convex}] \label{lem:KM-positions}
Let $K$ be a $2$-convex body with constant $\alpha$ and volume 1.
\begin{enumerate}
\item
If $K$ is in minimal mean-width position then:
\[
M^*(K) \leq C \sqrt{n} \min(\frac{1}{\sqrt{\alpha}},\log(1+n)) ~.
\]
\item
In fact, the same estimate on $M^*(K)$ is valid in John's
position.
\end{enumerate}
\end{lem}

The latter easily generalizes to the case of general $p$-convex
bodies. We sketch the argument for the following lemma (see
\cite{LT-Book-II} for definitions):

\begin{lem} \label{lem:general-positions}
Let $K$ be a $p$-convex body with constant $\alpha$ and volume 1.
If $K$ is in minimal mean-width position then:
\[
M^*(K) \leq \sqrt{n} \min(f(p,\alpha),C \log(1+n)) ~,
\]
where $f$ is a function depending solely on $p$ and $\alpha$.
\end{lem}

\begin{proof}[Sketch of proof]
Recall that by the classical result of Figiel and Tomczak-Jaegermann
on the $l$-position (\cite{l-position}), we have that in the minimal
mean-width position, $M^*(K) \leq C \sqrt{n} \norm{Rad(X_K)}$ for a
convex body $K$ of volume 1, where $\norm{Rad(X_K)}$ denotes the
norm of the Rademacher projection on $L_2(X_K)$ (see e.g.
\cite{Milman-Schechtman-Book} for definitions). Since $K$ is
$p$-convex with constant $\alpha$, it is classical
(\cite[Proposition 1.e.2]{LT-Book-II}) that $K^\circ$ is $q$-smooth
($q = p^*$) with constant $\beta(\alpha,p)$, and therefore
(\cite[Theorem A.7]{BL-Book}) has type-$q$, with $T_q(X_K^*)$
depending only on $p$ and $\alpha$. Pisier showed in
\cite{Pisier-Type-Implies-K-Convex} that $\norm{Rad(X)} =
\norm{Rad(X^*)}$ may be bounded from above by an (explicit) function
of $T_q(X^*)$ when $q>1$, which shows that $M^*(K) \leq \sqrt{n}
f(p,\alpha)$. By another important result of Pisier (e.g.
\cite{Milman-Schechtman-Book}), for an $n$-dimensional Banach space
$X$ one always has $\norm{Rad(X)} \leq C \log(1+n)$, showing that
$M^*(K) \leq \sqrt{n} C \log(1+n)$.
\end{proof}

Combining Lemmas \ref{lem:p-convex-inscribed-ball} and
\ref{lem:general-positions} with Proposition
\ref{prop:variance-concentration}, we get a concentration result for
$p$-convex bodies with $2 \leq p <4$. The concentration will be for
$T(K)$, the position which is ``half-way" (in the geometric mean
sense) between the isotropic position $K$ and the minimal mean-width
position $T^*T(K)$.

\begin{thm} \label{thm:p-convex-variance}
Let $K$ be a $p$-convex body with constant $\alpha$ for $2\leq p
<4$. Assume that $K$ is isotropic and of volume 1, and set
$q=p^*$. Then there exists a position $T(K)$ with $T \in SL(n)$,
such that:
\begin{enumerate}
\item
\[
\norm{T}_{op} \leq C \frac{n^{\frac{1}{2q}}}{\alpha^{\frac{1}{2p}}
L_K^{\frac{1}{2}}} \min(f(p,\alpha),\log(1+n))^{\frac{1}{2}} ~.
\]
\item
\begin{equation} \label{eq:Var-estimate}
(\Var^q_{T(K)}(|x|^2))^{1/q} \leq C n^{\frac{1}{4} + \frac{1}{q}}
\alpha^{-\frac{1}{p}} L_K \min(f(p,\alpha),\log(1+n)) ~.
\end{equation}
\item
Set $\rho^2 = \int_{T(K)} |x|^2 dx / n$. Then:
\[
\VolSet{x \in T(K) ; \abs{\frac{|x|}{\sqrt{n}} - \rho} \geq t \rho}
\leq 2 \exp\brac{ - \frac{ c L_K^{1/2} \alpha^{\frac{1}{2p}}
n^{\frac{3}{8} - \frac{1}{2q}} t^{\frac{1}{2}}}{
\min(f(p,\alpha),\log(1+n))^{1/2} } } ~.
\]
\end{enumerate}
\end{thm}
\begin{proof}
Since the isotropic and the minimal mean-width positions are defined
up to orthogonal rotations, we may find a positive definite $T \in
SL(n)$ so that $T^*T(K)$ is in minimal mean-width position, which by
Lemma \ref{lem:general-positions} and Proposition
\ref{prop:variance-concentration} gives (2). Since $diam(T^*T(K))
\leq C \sqrt{n} M^*(T^*T(K))$, we also have:
\begin{equation} \label{eq:John-Will-Improve}
T^*T(K) \subset C n \min(f(p,\alpha),\log(1+n)) D_n ~.
\end{equation}
By Lemma \ref{lem:p-convex-inscribed-ball}, this means that:
\[
\norm{T^* T}_{op} \leq C
\frac{n^{1-\frac{1}{p}}}{\alpha^{\frac{1}{p}} L_K}
\min(f(p,\alpha),\log(1+n)) ~,
\]
which gives (1).
To deduce (3), we use the results of Bobkov \cite{BobkovPolynomials}
on the growth of $L_r$ norms of polynomials. Note that the function
$g(x) = \abs{x}^2 - n \rho^2$ is a polynomial of degree 2, so by
\cite[Theorem 1]{BobkovPolynomials} there exists a universal
constant $C>0$ such that:
\begin{equation} \label{eq:Bobkov}
E_{T(K)}\brac{\exp\brac{\frac{|g|^{1/2}}{C E_{T(K)}(|g|^{1/2})}}}
\leq 2 ~.
\end{equation}
Since $E_{T(K)}(|g|^{1/2}) \leq E_{T(K)}(|g|^q)^{\frac{1}{2q}} =
\Var_{T(K)}^q(\abs{x}^2)^{\frac{1}{2q}}$, using the Chebyshev-Markov
inequality, (\ref{eq:Bobkov}) and (\ref{eq:Var-estimate}), yields:
\begin{eqnarray}
\nonumber & & \!\!\!\!\!\!\!\! \VolSet{x \in T(K) ;
\abs{\frac{|x|}{\sqrt{n}} - \rho} \geq t \rho} \leq \VolSet{x \in
T(K) ; \abs{\abs{x}^2 - n
\rho^2} \geq n \rho^2 t} \\
\nonumber & = &\VolSet{x \in T(K) ; |g(x)|^{1/2} \geq \sqrt{n t}
\rho} \leq 2 \exp\brac{-\frac{\sqrt{n t} \rho}{C
\Var_{T(K)}^q(|x|^2)^{1/2q}}} \\
\nonumber &\leq & 2 \exp\brac{ - \frac{ \rho \alpha^{\frac{1}{2p}}
n^{\frac{3}{8} - \frac{1}{2q}} t^{\frac{1}{2}}}{ C' L_K^{1/2}
\min(f(p,\alpha),\log(1+n))^{1/2} } } ~.
\end{eqnarray}
(3) immediately follows since always $\rho \geq L_K$ (e.g.
\cite{Milman-Pajor-LK}).
\end{proof}

\begin{rem} \hfill
\begin{enumerate}
\item
We see from (2) and (3) that we get a non-trivial concentration when
$q > \frac{4}{3}$, i.e. $p<4$; this is due to the extra
$n^{\frac{1}{4}}$ term in Proposition \ref{prop:variance-concentration}.
\item
For 2-convex bodies, we can slightly improve the estimate on
$\norm{T}_{op}$ by taking $T^*T(K)$ to be in John's position.
Indeed, by part (2) of Lemma \ref{lem:KM-positions}, we will have
the same estimate on $M^*(T^*T(K))$ as the one used in the proof of
Theorem \ref{thm:p-convex-variance}. The advantage of using John's
position is that $T^*T(K) \subset C n D_n$, improving the estimate
in (\ref{eq:John-Will-Improve}), which was used to derive the bound
on $\norm{T}_{op}$.
\end{enumerate}
\end{rem}

The advantage of this theorem over the previous concentration
results for $p$-convex bodies in \cite{ABP} or
\cite{KlartagEMilman-2-Convex} is three-fold. In \cite{ABP}, the
concentration was shown under certain assumptions on the diameter of
the bodies, which is not satisfied for some bodies (as shown in
\cite{KlartagEMilman-2-Convex} even for $p=2$). In
\cite{KlartagEMilman-2-Convex}, this restriction on the diameter was
removed for $p=2$, but the resulting concentration depended on an
implicit function $\lambda = \lambda(\alpha)$, which appeared in the
exponent of $n$. In Theorem \ref{thm:p-convex-variance} for the case
$2\leq p <4$, the restrictions on the diameter of the bodies are
removed, the dependence of the concentration on $\alpha$ is
explicit, and this dependence is not in the exponent in any of the
expressions.

\medskip

Since it is well know that $L_K \geq c$ (e.g.
\cite{Milman-Pajor-LK}), Theorem \ref{thm:p-convex-variance} yields
a concentration of the form (\ref{eq:Sodin-concentration}) required
to apply Theorem \ref{thm:Sodin} to $T(K)$.
It remains to evaluate $C_{iso}(T(K))$, taking into account the
remark after Theorem \ref{thm:extended-ABP}.
Since $\rho_{avg} \geq c \rho \geq c L_K$, using (1) from Theorem
\ref{thm:p-convex-variance} and $L_K \geq c$, we have:
\[
C_{iso}(T(K)) = \frac{\rho_{max}}{\rho_{avg}} \leq
\frac{\norm{T}_{op} L_K}{c L_K} \leq C
n^{\frac{1}{2q}}\alpha^{-\frac{1}{2p}}
\min(f(p,\alpha),\log(1+n))^{\frac{1}{2}} ~.
\]
Plugging everything into Theorem \ref{thm:Sodin}, we deduce Theorem
\ref{thm:2-4-convex}. Corollary \ref{cor:2-4-convex} is deduced by
using the estimates given in Lemma \ref{lem:sqlp}.

\begin{rem}
Sasha Sodin has brought to our attention that a recent result of S.
Bobkov (\cite{BobkovVar}) shows that all our concentration results
for uniformly convex bodies in fact imply isoperimetric inequalities
for these bodies (with respect to the Euclidean norm). In fact, in a
recent manuscript by Sodin and the author
\cite{EMilmanSodinIsoperimetryForULC}, we prove isoperimetric analogues
of the Gromov-Milman Theorem for uniformly convex bodies, which may
be used directly to obtain isoperimetric inequalities with respect
to the Euclidean norm, by employing the estimates in this note.
\end{rem}

\setlinespacing{0.85}
\bibliographystyle{amsplain}
\bibliography{../../../ConvexBib}

\def\cprime{$'$}
\providecommand{\bysame}{\leavevmode\hbox to3em{\hrulefill}\thinspace}
\providecommand{\MR}{\relax\ifhmode\unskip\space\fi MR }
\providecommand{\MRhref}[2]{%
  \href{http://www.ams.org/mathscinet-getitem?mr=#1}{#2}
}
\providecommand{\href}[2]{#2}
\begin{thebibliography}{10}

\bibitem{ABP}
M.~Anttila, K.~Ball, and I.~Perissinaki, \emph{The central limit problem for
  convex bodies}, Trans. Amer. Math. Soc. \textbf{355} (2003), no.~12,
  4723--4735.

\bibitem{ABV}
J.~Arias-de Reyna, K.~Ball, and R.~Villa, \emph{Concentration of the distance
  in finite-dimensional normed spaces}, Mathematika \textbf{45} (1998), no.~2,
  245--252.

\bibitem{BL-Book}
Y.~Benyamini and J.~Lindenstrauss, \emph{Geometric nonlinear functional
  analysis. {V}ol. 1}, American Mathematical Society Colloquium Publications,
  vol.~48, American Mathematical Society, Providence, RI, 2000.

\bibitem{BobkovVar}
S.~Bobkov, \emph{On isoperimetric constants for log-concave probability
  distributions}, to appear in GAFA Seminar Notes 2004-2005, 2005.

\bibitem{BobkovPolynomials}
S.~G. Bobkov, \emph{Remarks on the growth of {$L\sp p$}-norms of polynomials},
  Geometric aspects of functional analysis, Lecture Notes in Math., vol. 1745,
  Springer, Berlin, 2000, pp.~27--35.

\bibitem{BobkovLedoux}
S.~G. Bobkov and M.~Ledoux, \emph{From {B}runn-{M}inkowski to {B}rascamp-{L}ieb
  and to logarithmic {S}obolev inequalities}, Geom. Funct. Anal. \textbf{10}
  (2000), no.~5, 1028--1052.

\bibitem{BobkovZegarlinski}
S.~G. Bobkov and B.~Zegarlinski, \emph{Entropy bounds and isoperimetry}, Mem.
  Amer. Math. Soc. \textbf{176} (2005), no.~829, x+69.

\bibitem{Borell-logconcave}
Ch. Borell, \emph{Convex measures on locally convex spaces}, Ark. Mat.
  \textbf{12} (1974), 239--252.

\bibitem{Bourgain-LK}
J.~Bourgain, \emph{On the distribution of polynomials on high dimensional
  convex sets}, Geometric Aspects of Functional Analysis, Lecture Notes in
  Mathematics, vol. 1469, Springer-Verlag, 1991, pp.~127--137.

\bibitem{BrehmVoigtVogt}
U.~Brehm, P.~Hinow, H.~Vogt, and J.~Voigt, \emph{Moment inequalities and
  central limit properties of isotropic convex bodies}, Mathematische
  Zeitschrift \textbf{240} (2002), no.~1, 37--51.

\bibitem{BrehmVoigt}
U.~Brehm and J.~Voigt, \emph{Asymptotics of cross sections for convex bodies},
  Beitr\"age Algebra Geom. \textbf{41} (2000), no.~2, 437--454.

\bibitem{Davis-etal-Lemma}
W.~J. Davis, V.~D. Milman, and N.~Tomczak-Jaegermann, \emph{The distance
  between certain $n$-dimensional banach spaces}, Israel Journal of Mathematics
  \textbf{39} (1981), 1--15.

\bibitem{DiaconisFreedmanProjectionPursuit}
P.~Diaconis and D.~Freedman, \emph{Asymptotics of graphical projection
  pursuit}, Ann. Statist. \textbf{12} (1984), no.~3, 793--815.

\bibitem{DiaconisFreedmanDeFinetti}
\bysame, \emph{A dozen de {F}inetti-style results in search of a theory}, Ann.
  Inst. H. Poincar\'e Probab. Statist. \textbf{23} (1987), no.~2, suppl.,
  397--423.

\bibitem{l-position}
T.~Figiel and N.~Tomczak-Jaegermann, \emph{Projections onto {H}ilbertian
  subspaces of {B}anach spaces}, Israel J. Math. \textbf{33} (1979), no.~2,
  155--171.

\bibitem{FleuryGuedonPaourisCLP}
B.~Fleury, O.~Gu{\'e}don, and G.~Paouris, \emph{A stability result for mean
  width of $l_p$-centroid bodies}, Advances in Mathematics \textbf{214} (2007),
  no.~2, 865--877.

\bibitem{Gromov-Milman}
M.~Gromov and V.~D. Milman, \emph{Generalization of the spherical isoperimetric
  inequality to uniformly convex {B}anach spaces}, Compositio Math. \textbf{62}
  (1987), no.~3, 263--282.

\bibitem{LinnikBook}
I.~A. Ibragimov and Yu.~V. Linnik, \emph{Independent and stationary sequences
  of random variables}, Wolters-Noordhoff Publishing, Groningen, 1971, With a
  supplementary chapter by I. A. Ibragimov and V. V. Petrov, Translation from
  the Russian edited by J. F. C. Kingman.

\bibitem{KlartagUnconditionalVariance}
B.~Klartag, \emph{A {B}erry-{E}sseen type inequality for convex bodies with an
  unconditional basis}, to appear in Prob. Theo. Relat. Fields,
  arXiv:0705.0832, 2007.

\bibitem{KlartagCLP}
\bysame, \emph{A central limit theorem for convex sets}, Invent. Math.
  \textbf{168} (2007), 91--131.

\bibitem{KlartagCLPpolynomial}
\bysame, \emph{Power-law estimates for the central limit theorem for convex
  sets}, J. Funct. Anal. \textbf{245} (2007), 284--310.

\bibitem{KlartagEMilman-2-Convex}
B.~Klartag and E.~Milman, \emph{On volume distribution in $2$-convex bodies},
  to appear in Israel Journal of Mathematics, www.arxiv.org/math.FA/0604594,
  2006.

\bibitem{LT-Book-II}
J.~Lindenstrauss and L.~Tzafriri, \emph{Classical {B}anach spaces. {II}},
  Ergebnisse der Mathematik und ihrer Grenzgebiete [Results in Mathematics and
  Related Areas], vol.~97, Springer-Verlag, Berlin, 1979, Function spaces.

\bibitem{Meckes-GaussianMarginals}
E.S. Meckes and M.W. Meckes, \emph{The central limit problem for random vectors
  with symmetries}, Manuscript, www.arxiv.org/math.PR/0505618, 2005.

\bibitem{Meckes-GaussianMarginalsNew}
M.W. Meckes, \emph{Gaussian marginals of convex bodies with symmetries},
  Manuscript, www.arxiv.org/abs/math/0606073, 2007.

\bibitem{EMilman-DualMixedVolumes}
E.~Milman, \emph{Dual mixed volumes and the slicing problem}, Advances in
  Mathematics \textbf{207} (2006), no.~2, 566--598,
  www.arxiv.org/math.FA/0512207.

\bibitem{EMilmanGaussianMarginalsArxiv}
\bysame, \emph{On gaussian marginals of uniformly convex bodies}, Manuscript,
  www.arxiv.org/math.FA/0604595, 2006.

\bibitem{EMilmanSodinIsoperimetryForULC}
E.~Milman and S.~Sodin, \emph{An isoperimetric inequality for uniformly
  log-concave measures and uniformly convex bodies}, J. Funct. Anal.
  \textbf{254} (2008), no.~5, 1235--1268, www.arxiv.org/abs/math/0703857.

\bibitem{Milman-Pajor-LK}
V.~D. Milman and A.~Pajor, \emph{Isotropic position and interia ellipsoids and
  zonoids of the unit ball of a normed $n$-dimensional space}, Geometric
  Aspects of Functional Analysis, Lecture Notes in Mathematics, vol. 1376,
  Springer-Verlag, 1987-1988, pp.~64--104.

\bibitem{Milman-Schechtman-Book}
V.~D. Milman and G.~Schechtman, \emph{Asymptotic theory of finite-dimensional
  normed spaces}, Lecture Notes in Mathematics, vol. 1200, Springer-Verlag,
  Berlin, 1986, With an appendix by M. Gromov.

\bibitem{Paouris-IsotropicTail}
G.~Paouris, \emph{Concentration of mass on convex bodies}, to appear in GAFA,
  2006.

\bibitem{Pisier-Type-Implies-K-Convex}
G.~Pisier, \emph{Holomorphic semigroups and the geometry of {B}anach spaces},
  Ann. of Math. (2) \textbf{115} (1982), no.~2, 375--392.

\bibitem{SchechtmanZinn}
G.~Schechtman and J.~Zinn, \emph{On the volume of the intersection of two
  {$L\sp n\sb p$} balls}, Proc. Amer. Math. Soc. \textbf{110} (1990), no.~1,
  217--224.

\bibitem{SchmuckyUniformlyConvexBodies}
M.~Schmuckenschl{\"a}ger, \emph{A concentration of measure phenomenon on
  uniformly convex bodies}, Geometric aspects of functional analysis (Israel,
  1992--1994), Oper. Theory Adv. Appl., vol.~77, Birkh\"auser, Basel, 1995,
  pp.~275--287.

\bibitem{Sodin-GaussianMarginals}
S.~Sodin, \emph{Tail-sensitive gaussian asymptotics for marginals of
  concentrated measures in high dimension}, Manuscript,
  www.arxiv.org/math.MG/0501382, 2005.

\bibitem{SudakovMarginals}
V.~N. Sudakov, \emph{Typical distributions of linear functionals in
  finite-dimensional spaces of high dimension}, Dokl. Akad. Nauk SSSR
  \textbf{243} (1978), no.~6, 1402--1405.

\bibitem{Talagrand-Book}
M.~Talagrand, \emph{The generic chaining}, Springer-Verlag, Berlin, 2005, Upper
  and lower bounds of stochastic processes.

\bibitem{TJ-Schatten-Modulus}
N.~Tomczak-Jaegermann, \emph{The moduli of smoothness and convexity and the
  {R}ademacher averages of trace classes {$S\sb{p}(1\leq p<\infty )$}}, Studia
  Math. \textbf{50} (1974), 163--182.

\bibitem{TJ-Type-Cotype-fev-vectors}
\bysame, \emph{Computing {$2$}-summing norm with few vectors}, Ark. Mat.
  \textbf{17} (1979), no.~2, 273--277.

\bibitem{vonWeiz}
H.~von Weizs{\"a}cker, \emph{Sudakov's typical marginals, random linear
  functionals and a conditional central limit theorem}, Probab. Theory Related
  Fields \textbf{107} (1997), no.~3, 313--324.

\bibitem{CoordCorrelationForOrliczNorms}
J.~O. Wojtaszczyk, \emph{The square negative correlation property for
  generalized orlicz balls}, to appear in GAFA Seminar Notes, 2005.

\end{thebibliography}

\end{document}